\documentclass[12pt]{iopart}
\usepackage{cite}
\usepackage{iopams}
\usepackage[pdftex]{xcolor}
\def\numpartsApp{\addtocounter{equation}{1}%
     \setcounter{eqnval}{\value{equation}}%
     \setcounter{equation}{0}%
     \def\theequation{%
     \Alph{section}.\arabic{eqnval}{\it\alph{equation}}}}

\def\endnumpartsApp{\def\theequation{%
     \Alph{section}.\arabic{equation}}%
     \setcounter{equation}{\value{eqnval}}}

\begin{document} 

\date{2020 Feb 18}

\title{Green function solution of generalised boundary value problems}
\author{Vanik E. Mkrtchian$^*$%
\footnote[3]{vem@daad-alumni.de}
and
Carsten Henkel$^\sharp$%
}
\address{$^*$\ 
Institute for Physical Research, Armenian Academy of Sciences,
Ashtarak 0203, Republic of Armenia
\\
$^\sharp$\ 
University of Potsdam, Institute of Physics and Astronomy, 
Karl-Liebknecht-Str. 24/25, 14476 Potsdam, Germany}

\begin{abstract}
We construct an expression for the Green function of a differential operator 
satisfying nonlocal, homogeneous boundary conditions starting from the
fundamental solution of the differential operator. This also provides the
solution to the boundary value problem of an
inhomogeneous partial differential equation with inhomogeneous, nonlocal,
and linear boundary conditions. The construction generally applies for all
types of linear partial differential equations and linear boundary
conditions.
\end{abstract}


\submitto{\JPA}

\section*{Introduction}


There are a number of examples where the use of Feynman's path integral
approach has provided simplifications, new results or a better understanding
of already known facts. Because it has a completely different starting point
compared to canonical quantum field theory, the functional integral may
identify unknown aspects of a problem and increase understanding
by adding complementary considerations. One item on this list is our earlier
work on a simple one-dimensional boundary value problem involving the
Helmholtz operator with Robin boundary conditions at the endpoints of
an interval \cite{AdP}. We solved the problem using Feynmans' path
integral and introduced auxiliary fields to take into account the
boundary conditions. We came to an expression for the Green function of
the problem that is very natural to generalize to a broader class of
differential operators, dimensions, and boundary conditions. 

It is the purpose of this paper to present this generalization and a
few of its consequences. We construct a Green function for a differential
operator with homogeneous boundary conditions and use it to solve
in a simple, but general way a linear inhomogeneous equation involving that 
differential operator and constrained by inhomogeneous linear
boundary conditions. The latter are taken in the form of an boundary
integral and thus may be called non-local. 

Throughout the paper, we use the
traditional field theory language (the functional integral representation
is avoided). The issues of existence and uniqueness of the solutions
are left for future work, assuming for the moment that we deal with
well-posed problems in Hadamard's sense \cite{Hadam}. As potential applications
we have in mind the electrodynamic response of interfaces with spatially 
dispersive materials that play a key role in dispersion forces 
\cite{Casimir-physics}
and radiative heat transfer across sub-micron vacuum gaps
\cite{heat-transfer}. In that context, the Green function provides a
compact way to compute the elements of the energy-momentum tensor of the
field within a bounded region, for example between two plates. It contains
(multiply) reflected waves that arise from the boundary conditions
considered on the plates. We show here how to connect the reflection
amplitudes in a simple and general way to the imposed boundary conditions.

\section{Problem formulation}

We consider a differential operator $\mathfrak{L}$ defined on $\mathbb{R}^{n}$,
a domain $\Omega \subset \mathbb{R}^{n}$ 
and want to construct a Green function $G_{x'}: \Omega \to \mathbb{C}$ 
that solves
\begin{equation}
\mathfrak{L} \, G_{x'} = \delta_{x'}  
\label{Prel.1}
\end{equation}%
for a source point $x'$ in the interior of $\Omega$. $\delta_{x'}$ is
the unit-mass Dirac distribution supported at $x'$ \cite{Kanwal}.
We shall often write $G_{x'}( x ) = G( x, x' )$.

The Green function is constrained by a finite set of 
homogeneous, nonlocal boundary conditions ($j = 1, \ldots m$)
\begin{equation}
\forall {{\bar{x} \in \partial\Omega}}:
\quad
\int_{\partial\Omega}\kern-1.3ex{\rm d}S( \bar{x}_1 ) \, 
b_j( \bar{x}, \bar{x}_{1} )
G( \bar{x}_{1}, {x'} )
= 0
\label{Prel.2}
\end{equation}
In this expression, ${\rm d}S( \bar{x}_1 )$ is the Lebesgue measure on the
boundary domain $\partial\Omega \ni \bar{x}_1$, 
and the $b_j: \partial\Omega \times \partial\Omega \to \mathbb{C}$ 
are well-behaved integral kernels.
Generally speaking, each kernel $b_{j}$ is defined on its own sub-domain
$\partial \Omega_{j}$ of the boundary surface $\partial\Omega
= \partial \Omega_{1} \cup \partial \Omega_{2} \cup \ldots 
\cup \partial \Omega_{m}$. The surface integral~(\ref{Prel.2}) can then 
be restricted to the sub-domain $\partial \Omega_{j}$.

In other words, the set of kernels $\{ b_j | j = 1, \ldots m\}$ 
defines
a linear map (integral operator) $\hat{B}$ from the space of complex functions
defined 
on the boundary domain $\partial\Omega$ to ``spinor-valued'' functions
$\hat{\varphi}: \partial\Omega \to \mathbb{C}^{m}$.
The restriction $\bar G_{x'}: \partial\Omega \to \mathbb{C},
\bar G_{x'}( \bar{x} ) = G_{x'}( \bar{x} )$
of the Green function $G_{x'}$ 
is in the kernel of this map, $\hat{B} \bar G_{x'} = \hat{0}$.
In the following, we use the word ``boundary function'' for a function
defined on $\partial\Omega$.

We also need the adjoint $\hat{B}^{\dagger}$ of the map $\hat{B}$: it
maps spinor-valued boundary functions $\hat{\varphi}$
to scalar functions and is defined by the boundary integral
\begin{equation}
(\hat{B}^{\dagger}\hat{\varphi})( \bar{x} ) = 
\sum_j
\int_{\partial\Omega}\kern-1.3ex{\rm d}S( \bar{x}_1 ) \,
b_j^*( \bar{x}_1, \bar{x} ) 
\varphi_j( \bar{x}_1 )
\label{Prel.3}
\end{equation}
where $b_j^*( \bar{x}_1, \bar{x} ) = 
[b_j( \bar{x}_1, \bar{x} )]^*$ is the complex conjugate.

We assume that the fundamental solution $E$ of the differential operator
$\mathfrak{L}$ \cite{Stak}, i.e., the solution to 
\begin{equation}
\mathfrak{L} \,  E_{x'} = \delta_{x'}  
\label{Prel.4}
\end{equation}%
in $\mathbb{R}^{n}$ is known: $E_{x'}( x ) = E( x, x' )$.
Note that
$E_{x'}$ is actually a member of a class of solutions because it is not
constrained by boundary conditions.

Finally, we consider apart from the original boundary value 
problem~(\ref{Prel.1}, \ref{Prel.2}) also its adjoint version: the
corresponding Green function is denoted $G^a_{x'}$ 
and solves
\begin{equation}
\mathfrak{L}^a\, G^a_{x'} = \delta_{x'}  
\label{Prel.6.a}
\end{equation}%
in the domain $\Omega$. In terms of the natural sesquilinear 
form on $\Omega$, this means
\begin{equation}
\langle \mathfrak{L} u ,\, G^a_{x'} \rangle 
= \int_{\Omega}\!{\rm d}V( x ) [(\mathfrak{L} u)(x)]^* G^a(x, x')
= u^*(x')
\label{Prel.6.a.x}
\end{equation}
for all $x' \in \Omega$ and functions $u: \Omega \to \mathbb{C}$ in the
domain of $\mathfrak{L}$.
Here, the volume integral involves
the familiar Lebesgue measure ${\rm d}V(x)$ in the domain $\Omega$.
The boundary condition for the adjoint problem is
\begin{equation}
\hat{B}^a \bar{G}^a_{x'} = 0
\label{Prel.6.b}
\end{equation}
where $\hat{B}^a$ has also $m$ components with kernels $
b^a_j(\bar{x}, \bar{x}_1)$, as in (\ref{Prel.2}).
We shall write $E^a_{x'}$ for the fundamental solution to
(\ref{Prel.6.a}) in $\mathbb{R}^n$.

It is a well-known result of the theory of linear differential equations 
(see for example Refs.\,\cite{Stak}, \cite{Greenberg}) that the Green function
$G( x, x' )$ 
and its adjoint counterpart $G^{a}( x, x' )$ are linked by 
\begin{equation}
G^{a}( x, x' ) = G^{*}( x', x )
\,.
\label{Prel.7}
\end{equation}

\section{Boundary integral representation of the Green function}

We claim that the Green function $G_{x'}$ to the 
problem~(\ref{Prel.1}, \ref{Prel.2})
can be expressed via the fundamental solution $E$ in the following
way
\refstepcounter{equation}\label{GF.1}\addtocounter{equation}{-1}%
\numparts
\begin{equation}
\forall{x \in \Omega}: \quad
G( x, {x'} ) = E( x, x' )
- 
\int_{\partial\Omega}\kern-1.3ex{\rm d}S( \bar{x} ) \, 
E( x, \bar{x} )
J^{(1)}_{x'}( \bar{x} ) 
\label{GF.1.0}
\end{equation}
where the integral runs over the boundary $\partial\Omega$ and involves
the boundary function $J^{(1)}_{x'}: \partial\Omega \to \mathbb{C}$. 
The latter is constructed by applying a sequence of linear maps 
\begin{equation}
J^{(1)}_{x'} = 
\hat{B}_1^{\dagger}
\hat{\bf g}^{-1}
\hat{B}
\bar{E}_{x'}
\label{GF.1.b}
\end{equation}
to the restriction $\bar{E}_{x'}$ of the fundamental solution 
$E_{x'}$ to the boundary $\partial\Omega$.

The key element
in the construction~(\ref{GF.1}) is the inverse $\hat{\bf g}^{-1}$ 
of the $m\times m$ matrix operator $\hat{\bf g}$, 
a linear map between spinor-valued functions on the boundary.
The matrix elements of $\hat{\bf g}$ are given by the double
boundary integrals ($1 \le i,j \le m$)
\begin{equation}
g_{ij}\left( \bar{x}, \bar{x}'\right) 
=
\int_{\partial\Omega \times \partial\Omega}\kern-2ex
{\rm d}S( \bar{x}_1 ) {\rm d}S( \bar{x}_2 ) \,
b_i( \bar{x}, \bar{x}_{1} ) 
E( \bar{x}_{1}, \bar{x}_2 ) 
b_{1j}^*( \bar{x}', \bar{x}_2 )
\label{GF.2}
\end{equation}
\endnumparts%
where the kernels $b_{1j}$ are the components of the boundary operator
$\hat{B}_1$.

\paragraph{Proof.}
We apply the differential operator $\mathfrak{L}$ to both sides of
(\ref{GF.1.0}).
The fundamental solution in the first term gives 
the Dirac distribution $\delta_{x'}$. In the second term, we pull
$\mathfrak{L}$ under the boundary integral and get zero because $x'$ is 
in the interior of $\Omega$, while $\bar{x}_1$ is not.

The boundary condition~(\ref{Prel.2}) follows from the 
construction~(\ref{GF.2}) of the matrix operator $\hat{\bf g}$:
the application of $\hat{B}$ to the second term in (\ref{GF.1.0})
generates under the integral the operator product
$\hat{B} \bar{E} \hat{B}_1^\dagger = \hat{\bf g}$ applied to the spinor 
function $\hat{\varphi} = \hat{\bf g}^{-1} \hat{B} \bar{E}_{x'}$.
The boundary integral thus reduces to $\hat{B} \bar{E}_{x'}$ which
cancels the first term. 

Obviously the solution $G_{x'}$ exists if the kernel~(\ref{GF.2}) defines
an invertible map $\hat{\bf g}$ \cite{Stak}. 
We assume that this holds if we deal with a 
well-posed problem \cite{Hadam}.
\hspace*{\fill} $\square$

\bigskip\

The above demonstration is still lacking an explicit expression for the
operator $\hat{B}_1$. To find it, we first write down a similar solution
to the adjoint problem~(\ref{Prel.6.a}, \ref{Prel.6.b}):
\numparts
\begin{eqnarray}
G^a( x, {x'} ) &= E^a( x, x' )
- 
\int_{\partial\Omega}\kern-1.3ex
{\rm d}S( \bar{x} ) \, 
E^a( x, \bar{x} )
J^{(2)}_{x'}( \bar{x} ) 
\\
J^{(2)}_{x'}( \bar{x} ) &=
\hat{B}_2^{\dagger}
(\hat{\bf h})^{-1}
\hat{B}^a
\bar{E}^a_{x'}
\label{AGF.1.a}
\\
h_{ij}( \bar{x}, \bar{x}' ) &=
\int_{\partial\Omega \times \partial\Omega}\kern-2ex
{\rm d}S( \bar{x}_1 ) {\rm d}S( \bar{x}_2 ) \,
b^a_i( \bar{x}, \bar{x}_{1} ) 
E^a( \bar{x}_{1}, \bar{x}_2 ) 
b_{2j}^*( \bar{x}', \bar{x}_2 )
\label{AGF.1.b}
\end{eqnarray}
\endnumparts
where $\hat{B}_2$ is another boundary operator.
The proof that this solves the adjoint problem follows the same lines
as above. (It is useful to note the concise form
of (\ref{AGF.1.b}): $\hat{\bf h} = \hat{B}^a \bar{E} \hat{B}_2^\dagger$.)

Because of the ambiguity of the unbounded fundamental solutions,
we may always assume that they satisfy the identity~(\ref{Prel.7}), i.e.
\begin{equation}
E^a( x, x' ) = E^*( x', x )
\label{SCC}
\end{equation}
By requiring that the matrix elements
$h_{ij}$ [Eq.(\ref{AGF.1.b})]
and $g^*_{ji}$ [Eq.(\ref{GF.2})] coincide, 
we find component-wise the identification
\begin{equation}
\hat{B}_1 = \hat{B}^{a},
\qquad
\hat{B}_{2} = \hat{B}
\label{SCC.2}
\end{equation}
We insert this \emph{Ansatz} 
into $G^*( x', x)$
from the complex conjugate of (\ref{GF.1.0}) 
and switch the boundary integral over
$E^*( x', \bar{x} )
J^{(1)*}_{x}( \bar{x} )$
written there with the application of
$\hat{B} \bar{E}_{x}$ in (\ref{GF.1.b}).
It is then easy to check that one gets the identity
\begin{equation}
\hat{B}_2^\dagger (\hat{\bf h})^{-1} \hat{B}^a 
= \hat{B}^\dagger \hat{\bf g}^{-1\dagger} \hat{B}_1
= \hat{B}_1^{\dagger*} \hat{\bf g}^{-1*} \hat{B}^*
\label{SCC.1}
\end{equation}
The last step in~(\ref{SCC.1}) is based on the observation that
in the position representation, the two-variable kernel of
this expression is a scalar function. Therefore, one may
take the formal transpose of the operator product.
This yields (\ref{Prel.7}) between $G$ and $G^a$.

When the solution~(\ref{SCC.2}) is inserted into~(\ref{GF.1}), 
we find the explicit expressions
\refstepcounter{equation}\label{OGF.2}
\addtocounter{equation}{-1}
\numparts
\begin{eqnarray}
G( x, {x'} ) &= E( x, x' )
- 
\int_{\partial\Omega}\kern-1.3ex{\rm d}S( \bar{x} ) \, 
E( x, \bar{x} )
J^{(1)}_{x'}( \bar{x} ) 
\label{OGF.2.a}
\\
J^{(1)}_{x'} &= 
\hat{B}^{a\dagger}
\hat{\bf g}^{-1}
\hat{B}
\bar{E}_{x'}
\\
g_{ij}\left( \bar{x}, \bar{x}'\right) 
&=
\int_{\partial\Omega \times \partial\Omega}\kern-2ex
{\rm d}S( \bar{x}_1 ) {\rm d}S( \bar{x}_2 ) \,
b_i( \bar{x}, \bar{x}_{1} ) 
E( \bar{x}_{1}, \bar{x}_2 ) 
b^{a*}_{j}\left( \bar{x}', \bar{x}_2\right)
\label{OGF.2.b}
\end{eqnarray}
\endnumparts
\paragraph{Remark.}
It follows from this that the Green function also obeys
\begin{equation}
G(x, x') \, 
\overleftarrow{\mathfrak{L}}_{x'} 
= \delta( x - x' )
\label{eq:L-from-right}
\end{equation}
with differential operator and boundary condition acting from the right:
\begin{equation}
\int_{\partial\Omega}\kern-1.3ex{\rm d}S( \bar{x}' ) \, 
G( x, \bar{x}' )
\hat{B}^{a\dagger}_j( \bar{x}', \bar{x} )
= 0
\label{eq:BS-from-right}
\end{equation}

\paragraph{Example.}
As a simple application of the solution~(\ref{OGF.2}) to the boundary
value problem~(\ref{Prel.1}),
let us construct the Green function of a self-adjoint Dirichlet problem, i.e.%
\begin{equation}
\mathfrak{L}_{x}G_{D}( x, x') = \delta( x - x' )
\,,
\qquad x, x'\in \Omega   
\label{Dir.a}
\end{equation}%
\begin{equation}
G_{D}( \bar{x}, x' ) = G_{D}^{a}( \bar{x}, x' ) = 0,
\qquad \bar{x} \in \partial\Omega  
\label{Dir.b}
\end{equation}%
Let $E( x, x')$ be the fundamental solution to~(\ref{Dir.a}),
then we find from the definition (\ref{OGF.2.b}) that the matrix
operator $\hat{\bf g}$ reduces to a scalar kernel
\begin{equation}
g( \bar{x}, \bar{x}' ) = E( \bar{x}, \bar{x}' )   
\label{Dir.c}
\end{equation}%
The solution~(\ref{OGF.2}) reads in this case%
\begin{eqnarray}
G( x, x' ) &= E( x, x' )
\nonumber\\
&\qquad {} - 
\int_{\partial\Omega \times \partial\Omega}\kern-2ex
{\rm d}S( \bar{x} ) {\rm d}S( \bar{x}' ) \,
E( x, \bar{x} ) g^{-1}( \bar{x}, \bar{x}' )
E( \bar{x}', x' )   
\label{Dir.d}
\end{eqnarray}%
where $g^{-1}$ is the inverse of the boundary integral operator with
kernel~(\ref{Dir.c}) defined on the manifold $\partial\Omega$:
\begin{eqnarray}
\int_{\partial\Omega}\kern-1.3ex{\rm d}S(\bar{x}_{1}) \,
E( \bar{x}, \bar{x}_{1} ) 
g^{-1}(\bar{x}_{1}, \bar{x}')
= 
\nonumber\\
\int_{\partial\Omega}\kern-1.3ex{\rm d}S(\bar{x}_{1}) \, 
g^{-1}( \bar{x}, \bar{x}_{1} ) 
E( \bar{x}_{1}, \bar{x}' )
= \delta( \bar{x}, \bar{x}' ) 
\label{Dir.e}
\end{eqnarray}

\section{Boundary value problem}

The above construction provides an integral representation that solves
a boundary value problem for the inhomogeneous equation 
\begin{equation}
\mathfrak{L} \,  u = f
\label{BVP.1}
\end{equation}%
with a smooth function $f: \Omega \to \mathbb{C}$. The boundary conditions
can be cast in a fairly general form as a set of $m$ inhomogeneous 
integral equations for the restriction $\bar{u}$ of $u$ to the boundary:
\begin{equation}
\hat{B} \bar{u} = \hat{\Phi}
\label{BVP.2}
\end{equation}%
Here, the given ``spinor'' $\hat{\Phi}$ has the boundary
functions
$\Phi_{j}: \partial\Omega \to \mathbb{C}$ as its components. 
Explicitly, Eq.(\ref{BVP.2}) reads
\begin{equation}
\forall {{\bar{x} \in \partial\Omega}}:
\quad
\int_{\partial\Omega}\kern-1.3ex{\rm d}S( \bar{x}_1 ) \, 
b_j( \bar{x}, \bar{x}_{1} )
u( \bar{x}_{1} )
= \Phi_j( \bar{x} )
\qquad
\left(j = 1, \ldots m\right)
\label{BVP.2b}
\end{equation}
by direct analogy to~(\ref{Prel.2}).

We claim that when the boundary value 
problem defined by~(\ref{BVP.1}--\ref{BVP.2b}) is well-posed \`a 
la Hadamard, then its solution $u$ is given by the expression 
\refstepcounter{equation}\label{BVP.4}\addtocounter{equation}{-1}%
\numparts
\begin{equation}
u(x) =
\int_{\Omega}{\rm d}V(x_1) \, G( x, x_1 ) f( x_1 )
+ \int_{\partial\Omega}\kern-1.3ex{\rm d}S( \bar{x} ) \,
      E( x, \bar{x} ) J_{\hat{\Phi}}( \bar{x} )
\label{BVP.4.0}
\end{equation}
Under the boundary integral (second line),
the complex-valued boundary function $J_{\hat{\Phi}}$ 
depends linearly on the $\Phi_j$'s:
\begin{equation}
J_{\hat{\Phi}} = \hat{B}^{a\dagger} \hat{\bf g}^{-1} \hat{\Phi}
\label{eq:Jphi-boundary-source}
\end{equation}
\endnumparts%
in close analogy to~(\ref{OGF.2.a}).

\paragraph{Proof.}
When acting on both sides of~(\ref{BVP.4.0}) with 
$\mathfrak{L}$, 
we find $f$ because of~(\ref{Prel.1}). 
Applying the operator $\hat{B}$ on the boundary restriction $\bar{u}$,
we find $\hat{\Phi}$ because of~(\ref{Prel.2})
and the definition~(\ref{OGF.2.b}) of the boundary operators 
$\hat{\bf g}$ and $\hat{\bf g}^{-1}$.
\hspace*{\fill} $\square$

\section{Examples}

\subsection{Local boundary conditions}

For a Cauchy initial value problem or for boundary value problems
with Dirichlet, Neumann or Robin boundary conditions,
the operator $\hat{B}$ in~(\ref{Prel.3}) acts in a local way:
\begin{equation}
(\hat{B}\bar{G}_{x'})( \bar{x} ) 
= \hat{B}( \bar{x} ) \bar{G}_{x'}( \bar{x} ) 
= 0
\label{eq:local-BC}
\end{equation}
The expression~(\ref{OGF.2}) for the Green function then reads
\refstepcounter{equation}\label{Loc.2}\addtocounter{equation}{-1}
\numparts
\begin{eqnarray}
G( x, x' ) &= E( x, x' )
\nonumber\\
& \qquad {} - 
\int_{\partial\Omega \times \partial\Omega}\kern-2ex
{\rm d}S( \bar{x}_1 ) {\rm d}S( \bar{x}_2 ) \, 
E( x, \bar{x}_1 )
\hat{K}( \bar{x}_1, \bar{x}_2 ) 
\bar{E}( \bar{x}_2, x' )
\\
\hat{K}( \bar{x}_1, \bar{x}_2 ) &=
\hat{B}^{a\dagger}(\bar{x}_1)
\hat{\bf g}^{-1}(\bar{x}_1, \bar{x}_2)
\hat{B}(\bar{x}_2)
\end{eqnarray}
Here, the local operators appear to the left and right of 
the $m\times m$ matrix kernel $\hat{\bf g}^{-1}(\bar{x}, \bar{x}')$.
The latter is the inverse of the matrix operator $\hat{\bf g}$ defined
on the boundary manifold $\partial\Omega$ in terms of the matrix elements
\begin{equation}
\hat{g}_{ij}( \bar{x}, \bar{x}' ) =
\hat{b}_i( \bar{x} )
E( \bar{x}, \bar{x}' ) 
\hat{b}^{a*}_j( \bar{x}' )   \label{Loc.3}
\end{equation}
\endnumparts%
The solution of the boundary value problem~(\ref{BVP.1})\ with local boundary
conditions%
\begin{equation}
\forall \bar{x} \in \partial\Omega:
\quad
\hat{B}( \bar{x} ) u(\bar{x}) = \hat{\Phi}( \bar{x} ) 
\label{Loc.4}
\end{equation}%
is%
\refstepcounter{equation}\label{Loc.5s}\addtocounter{equation}{-1}
\numparts
\begin{eqnarray}
u(x) &= \int_{\Omega}{\rm d}V(x_1) \, G( x, x_1 ) f( x_1 )
+ \int_{\partial\Omega}\kern-1.3ex{\rm d}S( \bar{x} ) \,
      E( x, \bar{x} ) \hat{B}^{a\dagger}(\bar{x})
      \hat{J}'_{\hat{\Phi}}( \bar{x} )
\label{Loc.5}
\\
\hat{J}'_{\hat{\Phi}} &= \hat{\bf g}^{-1}\hat{\Phi}
\end{eqnarray}
\endnumparts

\subsection{One-dimensional boundary value problem}

In the case of ordinary differential equations, we deal with a one-dimensional
problem defined on the interval $\Omega = \left( a,b \right)$.
In this case, the boundary manifold 
is the set of two endpoints $\partial\Omega=\left\{ a, b\right\}$. 
We then have local boundary conditions that may be expanded
in the form
\begin{equation}
a_{j0}G_{x'}( a ) + a_{j1} \partial_{x}G_{x'}(a)
+ 
b_{j0}G_{x'}( b ) + b_{j1} \partial_{x}G_{x'}(b)
= 0
\label{ODE.1}
\end{equation}%
where $a_{j0},a_{j1},b_{j0}$ and $b_{j1}$ ($j = 1, \ldots m$)
are constants whose values encode
whether this linear combination
corresponds to the Dirichlet, Neumann oder Robin type.
Obviously, the object $\hat{\bf g}$ in~(\ref{Loc.3}) reduces in this case 
to a numerical matrix rather than a matrix-valued integral operator.

The action of the operator $\hat{b}_{j}$ from the left
on the fundamental solution $E$ in~(\ref{Loc.2}) 
is in this context to be understood as
\begin{eqnarray}
& \hat{b}_{j}(\bar{x})E( \bar{x}, x' ) \rightarrow
\nonumber\\
& \qquad
a_{j0}E( a, x' ) + a_{j1}\partial_{x}E( a, x' )
+ 
b_{j0}E( b, x' ) + b_{j1}\partial_{x}E( b, x' )
\label{ODE.2.L}
\end{eqnarray}
where in $\partial_{x}E$, the differentiation is with respect to the
first argument of $E$.  
The action of the operator $\hat{b}_{j}^{a*}$ from the right
is defined as
\begin{eqnarray}
& E( x, \bar{x}' ) \hat{b}_{j}^{a*}(\bar{x}') \rightarrow 
\nonumber\\
& \qquad
a_{j0}^{a*}E( x, a ) + a_{j1}^{a*}\partial_{x'}E( x, a )
+
b_{j0}^{a*}E( x, b ) + b_{j1}^{a*}\partial_{x'}E( x, b )
\label{ODE.2.R}
\end{eqnarray}
where $a_{j0}^{a},a_{j1}^{a},b_{j0}^{a}$ and $b_{j1}^{a}$ are constants
appearing in the adjoint boundary conditions, and the derivative
$\partial_{x'}E$ is with respect to the second argument of $E$.

\section{Concluding remarks}

In a previous paper~\cite{AdP}, we constructed the electromagnetic
Green function
in a bounded domain subject to (nonlocal) boundary conditions at the interface
between spatially dispersive media.
We have shown here that this result can be generalized to give a
boundary integral representation 
for the Green function related to a broad class of linear partial differential
equations with linear homogeneous and nonlocal boundary conditions.
This Green function provides the solution to a boundary value problem 
for linear, inhomogeneous partial differential equations subject to nonlocal,
inhomogeneous conditions on the boundary manifold.

\appendix

\section{}

In this Appendix, we use the block matrix inversion formula \cite{Zhang} 
\begin{eqnarray}
\left( 
\begin{array}{cc}
\mathbf{A} & \mathbf{B} \\ 
\mathbf{C} & \mathbf{D}%
\end{array}%
\right) ^{-1} &=
\left[ 
\begin{array}{cc}
\mathbf{A}^{-1}+\mathbf{A}^{-1}\mathbf{BRCA}^{-1} \,,
	& 
	-\mathbf{A}^{-1}\mathbf{BR} 
\\ 
-\mathbf{RCA}^{-1} \,,
	& 
	\mathbf{R}%
\end{array}%
\right] 
\nonumber\\
\hspace*{5em}
\mathbf{R} &= \left( \mathbf{D-CA}^{-1}\mathbf{B}\right) ^{-1}
\label{Inverse}
\end{eqnarray}%
to find a simple expression for the inverse of the matrix operator
$\hat{\bf g}$ [Eq.(\ref{OGF.2.b})] and of the Green function 
[Eq.(\ref{OGF.2.a})]. The resulting expression for the Green function
is
\refstepcounter{equation}\label{GFB}\addtocounter{equation}{-1}
\numpartsApp
\begin{equation}
G(x, x') = E(x, x') - \sum\limits_{j=1}^{m}G_{j}(x, x')
\label{GFB.a}
\end{equation}%
with the recursive construction%
\begin{eqnarray}
G_{j}(x, x') &= \int_{\partial\Omega^{\times 4}} \kern-2ex
G^{( j-1 )}(x, \bar{y}_{1}) 
\, \hat{b}_{j}^{a*}( \bar{y}_{2}, \bar{y}_{1} ) 
\, g_{j}^{-1}( \bar{y}_{2}, \bar{y}_{2}')
\, \hat{b}_{j}( \bar{y}_{2}', \bar{y}_{1}' )
\nonumber\\
& \qquad \qquad {} \times
\, G^{( j-1 )}(\bar{y}_{1}', x')  
\label{GFB.b} \\
g_{j}( \bar{x}, \bar{x}' ) &= \int_{\partial\Omega^{\times 2}} \kern-2ex
\hat{b}_{j}( \bar{x}, \bar{y} )
\, G^{( j-1 )}(\bar{y}, \bar{y}')
\, \hat{b}_{j}^{a*}( \bar{x}', \bar{y}' )   
\label{GFB.c}\\
G^{( j-1 )} &= E - G_{1} - \ldots - G_{j-1}
\nonumber\\
G^{( 0 )} & \equiv E  
\nonumber
\end{eqnarray}%
\endnumpartsApp
In these expressions, an integral over $\partial\Omega$ 
with respect to every doubly appearing variable
$\bar{y}_1, \ldots \bar{y}'_2$ has to be performed.
The necessary condition is that the kernels $g_j$ can be inverted on
the boundary manifold $\partial\Omega$:
\begin{eqnarray}
\int_{\partial\Omega}\!{\rm d}S( \bar{x}_1 )
\, g_j( \bar{x}, \bar{x}_1 ) g_j^{-1}( \bar{x}_1, \bar{x}' )
=
\nonumber\\
\int_{\partial\Omega}\!{\rm d}S( \bar{x}_1 )
\, g_j^{-1}( \bar{x}, \bar{x}_1 ) g_j( \bar{x}_1, \bar{x}' ) 
=
\delta( \bar{x}, \bar{x}' )
\label{eq:invertible-gj}
\end{eqnarray}

\paragraph{Proof.}
To verify this result, we act in the following way: we start from $m=1$
where~(\ref{OGF.2}) reads%
\refstepcounter{equation}\label{m.1}\addtocounter{equation}{-1}
\numpartsApp
\begin{eqnarray}
G( x, x' ) & = E( x, x' ) - G_{1}( x, x' )   \label{m.1.a}
\\
G_{1}( x, x' ) &= \int_{\partial\Omega^{\times 4}} \kern-2ex
E( x, \bar{x}_{1} ) 
\, \hat{b}_{1}^{a*}(\bar{x}_{2}, \bar{x}_{1})
\, {g}_{11}^{-1}(\bar{x}_{2}, \bar{x}_{2}')
\, \hat{b}_{1}(\bar{x}_{2}', \bar{x}_{1}')
\, E( \bar{x}_{1}', x' )  
\label{m.1.b}
\\
{g}_{11}( \bar{x}, \bar{x}' ) &= 
\int_{\partial\Omega^{\times 2}} \kern-2ex
\hat{b}_{1}( \bar{x}, \bar{x}_{1}) 
\, E( \bar{x}_{1}, \bar{x}_{1}') 
\, \hat{b}_{1}^{a*}( \bar{x}', \bar{x}_{1}' ) 
\label{m.1.c}
\end{eqnarray}
\endnumpartsApp
We see that the expressions (\ref{m.1}) coincide with 
(\ref{GFB}) for $m=1$.

Introducing the shorthand based on~(\ref{OGF.2.b})
\begin{equation}
g_{ij} \equiv g_{ij}( \bar{x}, \bar{x}' ) =
\int_{\partial\Omega^{\times 2}} \kern-2ex
\hat{b}_{i}( \bar{x}, \bar{x}_{1}) 
\, E( \bar{x}_{1}, \bar{x}_{1}') 
\, \hat{b}_{j}^{a*}( \bar{x}', \bar{x}_{1}')   \label{g.matrix}
\end{equation}%
we find the inverse of the matrix operator $\hat{\bf g}$ for $m=2$ 
using the formula~(\ref{Inverse})%
\begin{equation}
\left( 
\begin{array}{cc}
g_{11} & \, g_{12} \\ 
g_{21} & \, g_{22}%
\end{array}%
\right) ^{-1}=\left[ 
\begin{array}{cc}
g_{11}^{-1} + g_{11}^{-1}g_{12} R\, g_{21}g_{11}^{-1} \,,
& \, -g_{11}^{-1}g_{12} R 
\\ 
- R\, g_{21}g_{11}^{-1} \,,
& R%
\end{array}%
\right]   
\label{m.2.a}
\end{equation}%
where 
$R = \left( g_{22} - g_{21} g_{11}^{-1} g_{12} \right)^{-1}
= g_{2}^{-1}$ [see Eq.(\ref{GFB.c})].
And then, inserting (\ref{m.2.a}) into expression~(\ref{OGF.2.a}) for the 
Green function, we find
\begin{equation}
G = E - G_{1} - G_{2}
\end{equation}%
which is coincident with (\ref{GFB.a}) for $m=2$ and where $G_{2}$ is
defined by (\ref{GFB.b}) for $j=2$.

In the case of $m=3$ we deal with $3\times 3$ $\hat{g}$ matrix which we
invert using (\ref{Inverse}) taking
\begin{equation}
\eqalign{%
\mathbf{A} \rightarrow \left( 
\begin{array}{cc}
g_{11} & \, g_{12} \\ 
g_{21} & \, g_{22}%
\end{array}%
\right) 
\qquad
& \mathbf{B} \rightarrow \left( 
\begin{array}{c}
g_{13} \\ 
g_{23}%
\end{array}%
\right) 
\cr
\mathbf{C} \rightarrow \left( 
\begin{array}{cc}
g_{31} & g_{32}%
\end{array}%
\right) 
\qquad
& \mathbf{D} \rightarrow g_{33}
}
\end{equation}
After some algebra, we arrive again at (\ref{GFB}) for $m=3$.
\hspace*{\fill} $\square$

\paragraph{Remark.}
And finally, let us consider the case of local boundary conditions 
where the operator $\hat{b}_{i}$\ acts in the domain $\partial\Omega_{i}$ 
of the boundary manifold $\partial\Omega$. 
Then, two integrals drop out in~(\ref{GFB.b}, \ref{GFB.c}),
and we find%
\numpartsApp
\begin{eqnarray}
G_{j}(x, x') &= 
\int_{\partial\Omega_{j}^{\times 2}}\kern-2ex%
G^{(j-1)}(x, \bar{y}_{1}) \,
\hat{b}_{j}^{a*}( \bar{y}_{1} ) \,
g_{j}^{-1}( \bar{y}_{1}, \bar{y}_{1}' ) \,
\hat{b}_{j}( \bar{y}_{1}' ) \,
G^{( j-1)}(\bar{y}_{1}', x')  
\label{GFBL.1}
\\
g_{j}( \bar{x}, \bar{x}' ) &= 
\hat{b}_{j}( \bar{x} ) \,
G^{( j-1)}(\bar{x}, \bar{x}') \,
\hat{b}_{j}^{a*}( \bar{x}' )   
\label{GFBL.2}
\end{eqnarray}
\endnumpartsApp

\ack

We thank Nikolai Tarkhanov for helpful remarks.

%

\end{document}